# A Warm Restart Strategy for Solving Sudoku by Sparse Optimization Methods

Yuchao Tang[1], Zhenggang Wu[2], Chuanxi Zhu[1]

1. Department of Mathematics, Nanchang University, Nanchang 330031, P.R. China
2. School of Data and Computer Science, Sun Yat-Sen University,Guangzhou 510006, P.R. China

**Abstract:** This paper is concerned with the popular Sudoku problem. We proposed a warm restart strategy for solving Sudoku puzzles, based on the sparse optimization technique. Furthermore, we defined a new difficulty level for Sudoku puzzles. The efficiency of the proposed method is tested using a dataset of Sudoku puzzles, and the numerical results show that the accurate recovery rate can be enhanced from 84%+ to 99%+ using the L1 sparse optimization method.

**Keywords:** Sudoku; Sparse optimization; L1 norm; Linear programming.

## 1. Introduction

Sudoku is a popular numbers game worldwide. The classical Sudoku puzzle is played on a $9 \times 9$ grid which is broken down into nine $3 \times 3$ blocks that do not overlap. A Sudoku puzzle usually comes with a partially filled grid. The objective is to fill the $9 \times 9$ grid with the digits 1 to 9, so that each column and row, and each of the nine $3 \times 3$ sub-grids, contain the digits 1 through 9 exactly once. There must be a unique solution to a standard Sudoku puzzle. Figure 1 shows an example of a Sudoku puzzle and its solution.

| | | | | | | | | |
|---|---|---|---|---|---|---|---|---|
| | | | | | | | | 1 |
| | | | | | | | 2 | |
| | | | | | 3 | | | |
| | | | | 4 | | 5 | | |
| | | 6 | | | | 3 | | |
| | | 7 | 8 | 1 | | | | |
| | 1 | | | 2 | | | | 4 |
| | 3 | | | | | | 7 | |
| 9 | 5 | | | | | | | |
| 2 | 7 | 3 | 9 | 6 | 4 | 8 | 5 | 1 |
| 4 | 6 | 9 | 1 | 5 | 8 | 7 | 2 | 3 |
| 1 | 8 | 5 | 2 | 7 | 3 | 4 | 6 | 9 |
| 8 | 2 | 1 | 3 | 4 | 6 | 5 | 9 | 7 |
| 5 | 4 | 6 | 7 | 9 | 2 | 3 | 1 | 8 |
| 3 | 9 | 7 | 8 | 1 | 5 | 2 | 4 | 6 |
| 7 | 1 | 8 | 5 | 2 | 9 | 6 | 3 | 4 |
| 6 | 3 | 2 | 4 | 8 | 1 | 9 | 7 | 5 |
| 9 | 5 | 4 | 6 | 3 | 7 | 1 | 8 | 2 |

Figure 1. A standard Sudoku puzzle (left) and its solution (right).

Sudoku is an interesting problem. It has attracted considerable attention from mathematician and computer scientists. A number of methods for solving Sudoku have been proposed, such as recursive backtracking [1], simulated annealing [2], integer programming [3], Sinkhorn balancing algorithms [4], sparse optimization method [5,6] and the alternating projection method [7]. These methods have their own merit and their limitations. In particular, Bartlett et al. [3] represented the Sudoku problem as a binary integer linear programming problem and, solved it by using the built-in MATLAB "bintprog" function. However, this approach is time consuming, especially for solving some difficult Sudoku puzzles. Babu et al. [5] first solved Sudoku puzzles based on the sparse optimization method, which showed very promising results in dealing with Sudoku puzzles. In brief, they introduced a method for transforming a Sudoku puzzle into a linear system. Because the resulting linear system is under-determined, there exist an infinite number of solutions. They proved that the sparsest solution of this linear equation is the solution of the Sudoku. Furthermore, they suggested using the $\ell_1$-norm minimization and weighted $\ell_1$-norm minimization models to approximate the sparsest solution. They tested both models on many Sudoku examples, of varying levels of difficulty. Most of these could be solved, but their methods also failed on some Sudoku puzzles. In 2014, McGuire et al. [8] proved that there is no possible Sudoku puzzles with fewer than 17 numbers already filled in. This answered a long-term open problem in Sudoku, regarding the smallest number of clues that a Sudoku puzzle can have. Therefore, we applied sparse optimization methods to a Sudoku dataset with each puzzle has 17 known clues, containing a total of 49,151 numbers. The sparse optimization method achieved almost an 84% precision. This motivated us to investigate whether there exists a way to improve these sparse optimization methods to obtain a higher precision. To achieve this goal, we developed some warm restart techniques for solving Sudoku puzzles based on sparse optimization methods. The idea is to delete repeated numbers in the Sudoku solution, and then solve it again. Our method is tested numerically on a dataset of Sudoku puzzles of various levels of difficulty, ranging from easy to very hard. The numerical results are very promising.

This paper is organized as follows. In the next section, we address the problem of encoding Sudoku puzzles into an under-determined linear system, and subsequently solving them using $\ell_1$-norm and weighted $\ell_1$-norm sparse optimization methods. In section 3, we propose an improved strategy for solving Sudoku puzzles based on the sparse optimization methods. These strategies could be viewed as restart techniques. Numerical experiments are presented in order to demonstrate the efficiency of our proposed methods in section 4. Finally, we present some conclusions.

## 2. Sparse optimization models and methods for solving Sudoku

In this section, we review how a Sudoku puzzle can be formulated as a linear system. Then, we propose two sparse models for solving it. Finally, two kinds of linear programming methods are proposed for solving the corresponding sparse optimization model.

### 2.1 Sudoku puzzles represented as linear systems

Babu et al. [5] first introduced a method for transforming a Sudoku puzzle into a linear system. An implicit approach was proposed in [3]. For the sake of completeness, we present the details here. To code a Sudoku puzzle as a linear system, we require binary variables (i.e., 0/1) to code the integer numbers 1 to 9.

Table 1. The integer numbers 1 to 9 are coded by a nine-dimensional binary vector.

| Integer numbers | Binary vector | Integer numbers | Binary vector | Integer numbers | Binary vector |
|---|---|---|---|---|---|
| 1 | (1,0,0,0,0,0,0,0,0) | 4 | (0,0,0,1,0,0,0,0,0) | 7 | (0,0,0,0,0,0,1,0,0) |
| 2 | (0,1,0,0,0,0,0,0,0) | 5 | (0,0,0,0,1,0,0,0,0) | 8 | (0,0,0,0,0,0,0,1,0) |
| 3 | (0,0,1,0,0,0,0,0,0) | 6 | (0,0,0,0,0,1,0,0,0) | 9 | (0,0,0,0,0,0,0,0,1) |

Each entry in a $9 \times 9$ Sudoku puzzle is associated to a nine dimensional variable shown in the Table 1. This results in 729 variables. Here, we let $x_{729\times 1}$ denote a solution of a Sudoku puzzle. Then this must satisfy the Sudoku constraints. For example, the first row should comprise all of the numbers 1,…,9, which can be expressed as

$$(\underbrace{I_{9\times 9} I_{9\times 9} \cdots I_{9\times 9}}_{9}\, 0_{9\times 648})x = 1_{9\times 1},$$

where $I_{9\times 9}$ denotes the $9\times 9$ identity matrix, $0_{9\times 648}$ denotes a matrix of size $9\times 648$ with all elements equal to zero, and $1_{9\times 1}$ denotes a column vector with all elements equal to one. Similarly, the constraint that the first column should contain all the digits 1,…,9 can be expressed as follows:

$$\underbrace{(I_{9\times 9}0_{9\times 72}I_{9\times 9}0_{9\times 72}\cdots I_{9\times 9}0_{9\times 72})}_{9}x = 1_{9\times 1}.$$

The constraint that the $3\times 3$ box in the top-left corner should contain all the digits 1,…,9 can be expressed as follows:

$$(I_{9\times 9}I_{9\times 9}I_{9\times 9}0_{9\times 54}I_{9\times 9}I_{9\times 9}I_{9\times 9}0_{9\times 54}I_{9\times 9}I_{9\times 9}I_{9\times 9}0_{9\times 540})x = 1_{9\times 1}.$$

The constraint that the first cell should be filled can be expressed as

$$(\underbrace{11\cdots 1}_{9}\underbrace{00\cdots 0}_{720})x = 1.$$

Finally, the clues can also be expressed using linear equality constraints. For example, in Figure 1, the clue in the second row and the eighth column with the value 2 can be expressed as

$$(\underbrace{00\cdots 0}_{144}\ 010000000\ \underbrace{00\cdots 0}_{576})x = 1.$$

By combining all these constraints, the linear equality constraints on $x$ can be expressed in a generic form as

$$Ax = \begin{bmatrix} A_{row} \\ A_{col} \\ A_{box} \\ A_{cell} \\ A_{clue} \end{bmatrix} x = b = \begin{bmatrix} 1 \\ 1 \\ \vdots \\ 1 \\ 1 \end{bmatrix}, \tag{1}$$

where $A_{row}$, $A_{col}$, $A_{box}$, $A_{cell}$ and $A_{clue}$ denote the matrices associated with the different constraints of Sudoku puzzles. For $9\times 9$ Sudoku puzzles, the size of $A$ is $(324+N_C)\times 729$, where $N_C$ denotes the number of clues. For example, for the Sudoku in Figure 1, the size of A is $341\times 729$, and hence the linear system of equations (1) is under-determined, and has infinitely many solutions. However, not every solution of (1) is a valid solution of the Sudoku puzzle. Babu et al. [5] proved if the Sudoku puzzle has a unique solution, then the sparsest solution of (1) is the solution for the Sudoku puzzle.

### 2.2 Sparse optimization models

In order to find the sparsest solution of the linear equation (1), we consider the $\ell_0$ minimization problem

$$(P_0)\ \min\ \|x\|_0 \quad \text{s.t. } Ax = b,$$

where $\|x\|_0$ represents the number of nonzero elements in $x$. In fact, this $\ell_0$ minimization problem is a basic problem of compressive sensing. Because the $\ell_0$ minimization problem is an NP-hard and non-convex problem, a popular approach is to replace $\ell_0$ by the $\ell_1$-norm. Then, the $\ell_0$ minimization problem becomes the following $\ell_1$-norm minimization problem:

$$(P_1)\ \min\ \|x\|_1 \quad \text{s.t. } Ax = b,$$

where $\|x\|_1 = \sum_{i=1}^{n}|x_i|$. The $\ell_1$-norm minimization problem is convex, and hence it has a unique solution. To enhance the sparsity of the solution searched by $(P_1)$, Candes et al. [9] proposed a weighted $\ell_1$-norm minimization problem:

$$(WP_1)\ \min\ \|Wx\|_1$$

$$\text{s.t. } Ax = b,$$

where $W = \mathrm{diag}(w_1, w_2, \cdots w_n)$ is a diagonal matrix with $w_k = \frac{1}{|x_k|^{i-1} + \epsilon}$, $0 < \epsilon < 1$, $k = 1,2,\cdots n$, and $i$ is the iteration number. The weighted matrix $W$ is obtained by solving the problem ($P_1$) from the previous iteration. Therefore, ($WP_1$) could be viewed as solving a series of $\ell_1$-norm optimization problems ($P_1$).

### 2.3 Sparse Optimization methods

The $\ell_1$-norm optimization problem ($P_1$) is equivalent to a linear programming problem, which can be efficiently solved using many well-known software, such as MATLAB or Lingo. We show that ($P_1$) can be solved using the following two linear programming methods. First,

$$\text{(LP1)} \ \min \ 1^T \begin{bmatrix} \hat{x} \\ \check{x} \end{bmatrix}$$

$$\text{s.t. } [A \quad -A] \begin{bmatrix} \hat{x} \\ \check{x} \end{bmatrix} = b$$

$$\begin{bmatrix} \hat{x} \\ \check{x} \end{bmatrix} \in C,$$

where $\hat{x}_{n\times 1}, \check{x}_{n\times 1}$ and $C = \{x | x \geq 0\}$. The solution of ($P_1$) is obtained by letting $x = \hat{x}_{n\times 1} - \check{x}_{n\times 1}$. If the solution of ($P_1$) is nonnegative, then we have the second linear programming problem:

$$\text{(LP2)} \ \min \ 1^T x$$

$$\text{s.t. } Ax = b$$

$$x \in C,$$

where the constraint $C$ is the same as (LP1).

Next, we show that the problem ($WP_1$) can be solved via the (LP1) and (LP2).

---

(Weighted LP1) Solving the weighted $\ell_1$-norm minimization problem ($WP_1$) by (LP1)

---

Input: $L = 10, \hat{x} = 0, \check{x} = 0, x_{ori} = \hat{x} - \check{x}, tol = 1 \times 10^{-10}$;

For i = 1:L

$W = \mathrm{diag}\left(\frac{1}{|x_{ori}| + \epsilon}\right)$;

Based on the method (LP1), $x_{new} = \hat{x} - \check{x}$ is obtained by solving the following linear programming problem:

$$\min \ W \begin{bmatrix} \hat{x} \\ \check{x} \end{bmatrix}, \text{ s.t. } [A \quad -A] \begin{bmatrix} \hat{x} \\ \check{x} \end{bmatrix} = b, \ \begin{bmatrix} \hat{x} \\ \check{x} \end{bmatrix} \in C;$$

if $\|x_{new} - x_{ori}\| < tol$

break;

else

$x_{ori} = x_{new}$;

End

End

Output: $x_{new}$

---

(Weighted LP2) Solving the weighted $\ell_1$-norm minimization problem ($WP_1$) by (LP2)

---

Input: $L = 10, x_{ori} = 0, tol = 1 \times 10^{-10}$;

For i = 1:L

$W = \mathrm{diag}\left(\frac{1}{|x_{ori}|+\epsilon}\right)$;

Based on the method (LP2), $x_{new}$ is obtained by solving the following linear programming problem:

min $Wx$, s.t. $Ax = b$, $x \in C$;

if $\|x_{new} - x_{ori}\| < tol$

break;

else

$x_{ori} = x_{new}$;

End

End

Output: $x_{new}$

The $729 \times 1$ vector $x$ is a stack of 81 $9 \times 1$ sub-vectors, one for each of the 81 cells in the Sudoku puzzle. Each $9 \times 1$ sub-vector is represented by a Boolean-type value, being zero except for a 1 in the position of the digit assigned to that cell. In practice, the $\ell_1$-norm minimization problem ($P_1$) and weighted $\ell_1$-norm minimization problem ($WP_1$) cannot always find a solution that takes a value of exactly 0/1. Thus, we need to transform the corresponding sub-vector.

For example, if we have a sub-vector $e = (0.1,0.11,0.3,0.4,0.22,0.211,0.113,0.122,0.33)$, then we have that

$$e_i = \begin{cases} 1, i = I(\max(e_i)), \\ 0, \text{otherwise}, \end{cases}$$

where $I(\max(e_i))$ denotes the position of the maximum number of in the vector $e$. Then, we obtain a new sub-vector $e = (0,0,0,1,0,0,0,0,0)$ which represents the integer number 4.

## 3. A warm restart strategy for solving Sudoku

If the solved Sudoku puzzle is wrong, this means that there is at least one row, column or $3 \times 3$ sub-square containing repeated numbers. Then, we propose an improved strategy in order to obtain the correct Sudoku solution. Our strategies can be divided into following steps:

**Step 1. (Second Solver)** Delete repeated numbers appearing in a row, column or $3 \times 3$ sub-square. Keep the remaining numbers, and view the results as a new Sudoku puzzle. If it is solved successfully, then quit. Otherwise, go to Step 2.

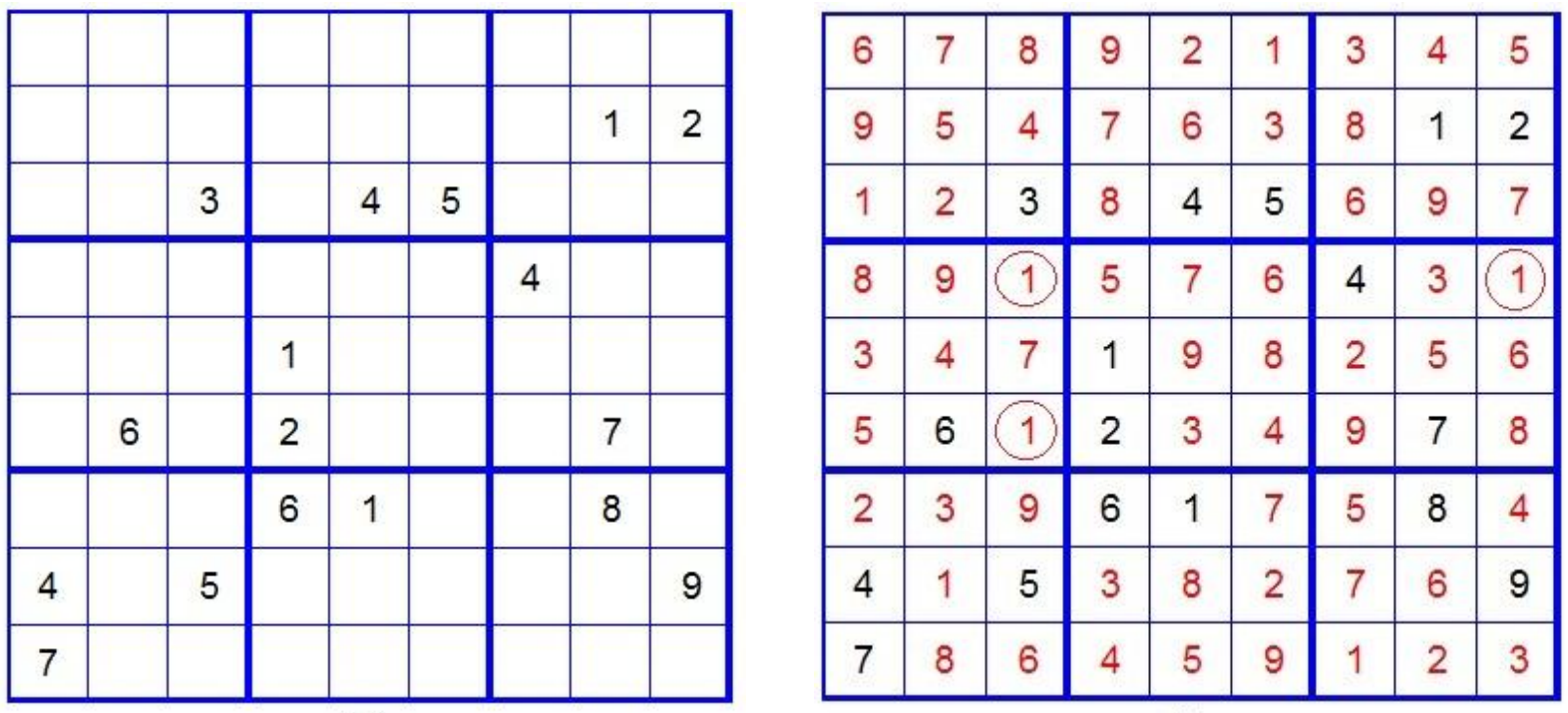

(a) (b)

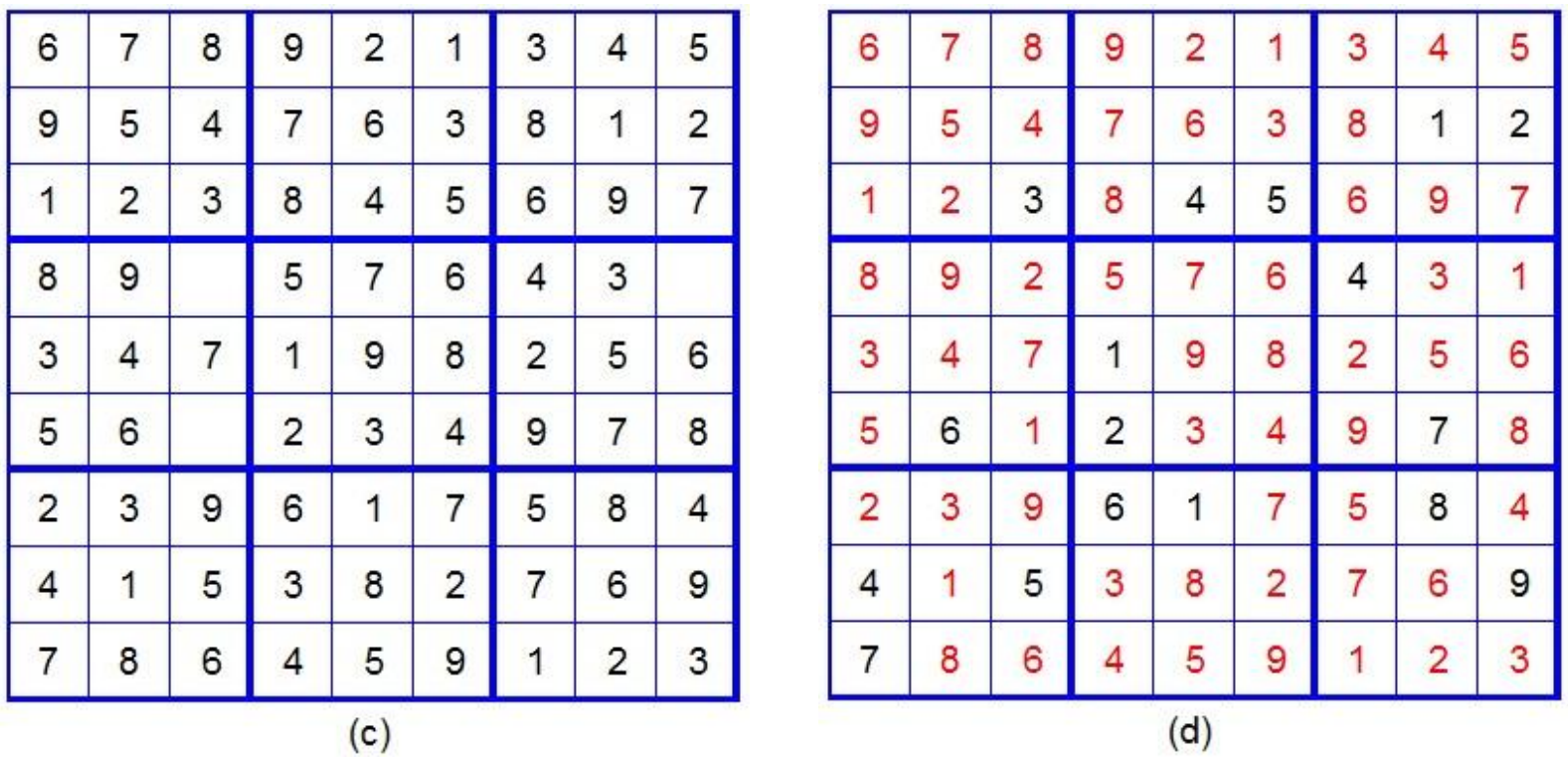

Figure 2. (a)The original Sudoku puzzle; (b) First solved solution, which is wrong: the wrong numbers are indicated by circles. (c) New Sudoku puzzle obtained by deleting the wrong numbers in (b); (d) Final correct solution.

**Step 2. (Third Solver)** Continue to delete repeated numbers appearing in a row, column or $3 \times 3$ sub-square. We can again obtain a new Sudoku puzzle, and then if it is solved successfully, stop. Otherwise, turn to Step 3.

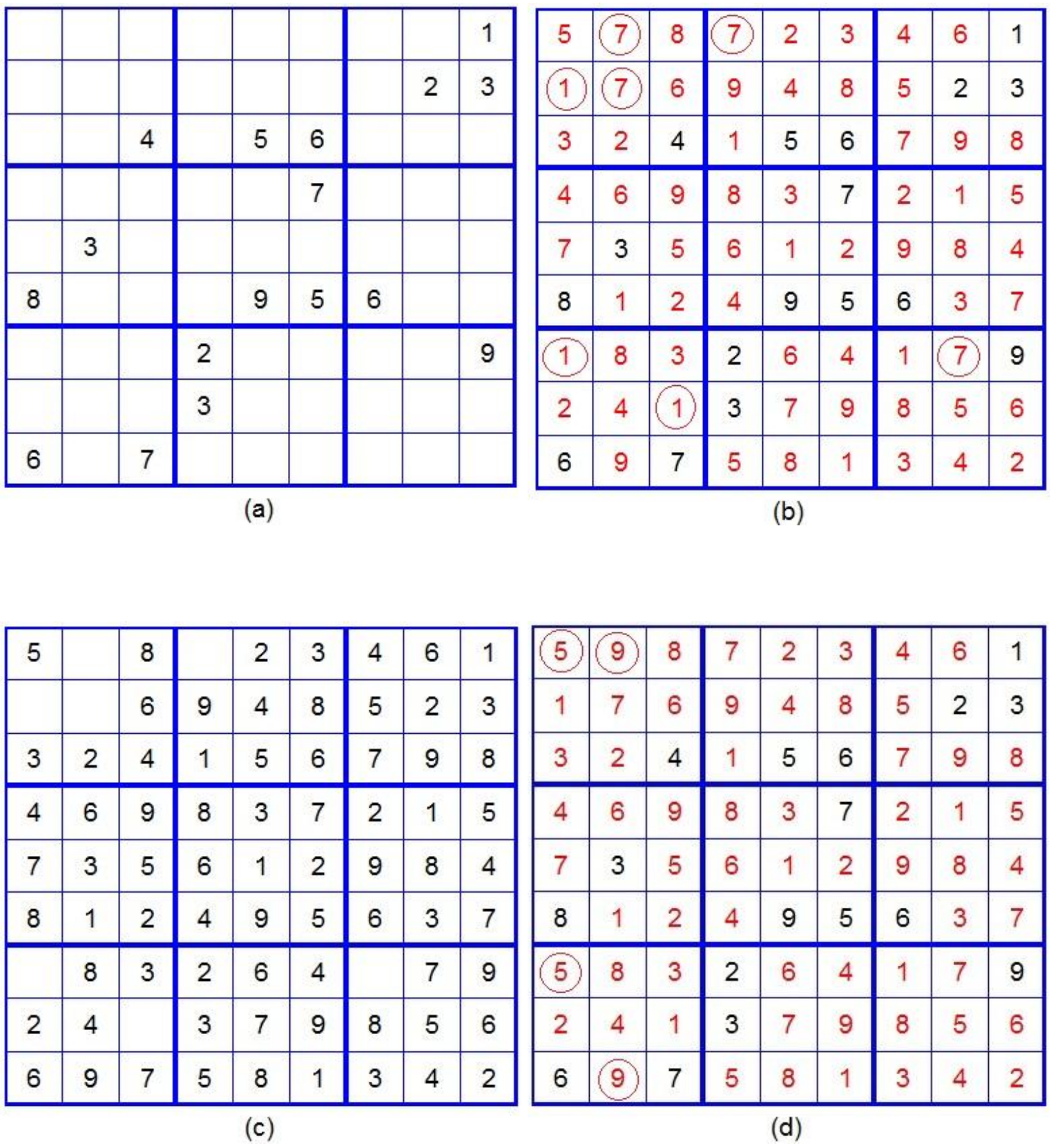

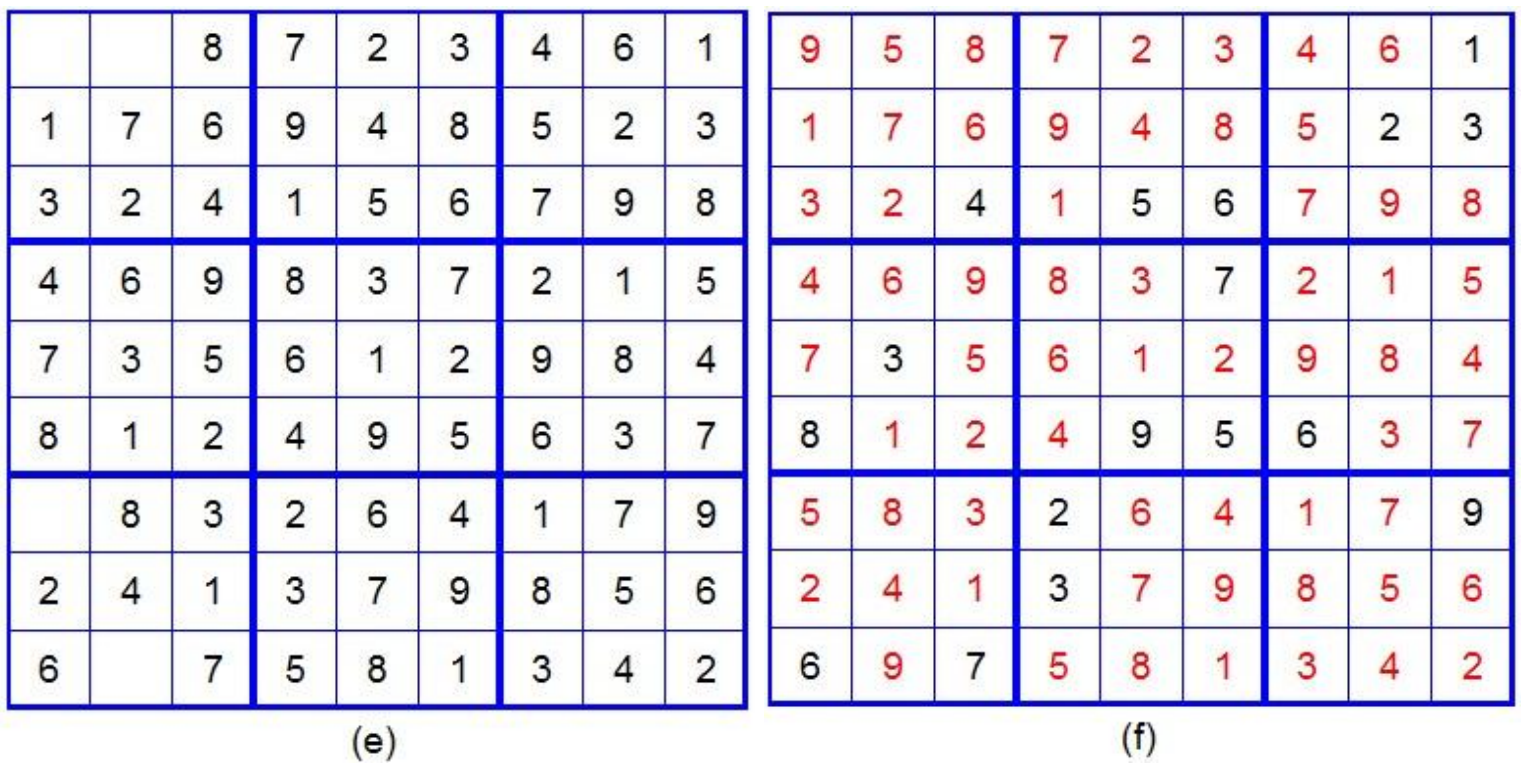

Figure 3. (a) The original Sudoku puzzle; (b) First solved solution, which is wrong; (c) New Sudoku puzzle obtained by deleting the wrong numbers in (b); (d) Second solved solution, which is still not correct; (e) New Sudoku puzzle obtained by deleting the wrong numbers in (d); (f) Final correct solution.

We found that if we repeat Step 2 again, there is no improvement. Thus, we turn to the third step.

**Step 3. (Successive increase Solver)** First, we delete repeated numbers appearing in a row, column or $3 \times 3$ sub-square, and remove the original Sudoku clue numbers. Then, select one number added to the original Sudoku puzzle, and hence obtain a newly defined Sudoku puzzle. If this Sudoku is solved correctly, quit. Otherwise, stop and return that this Sudoku is unsolvable.

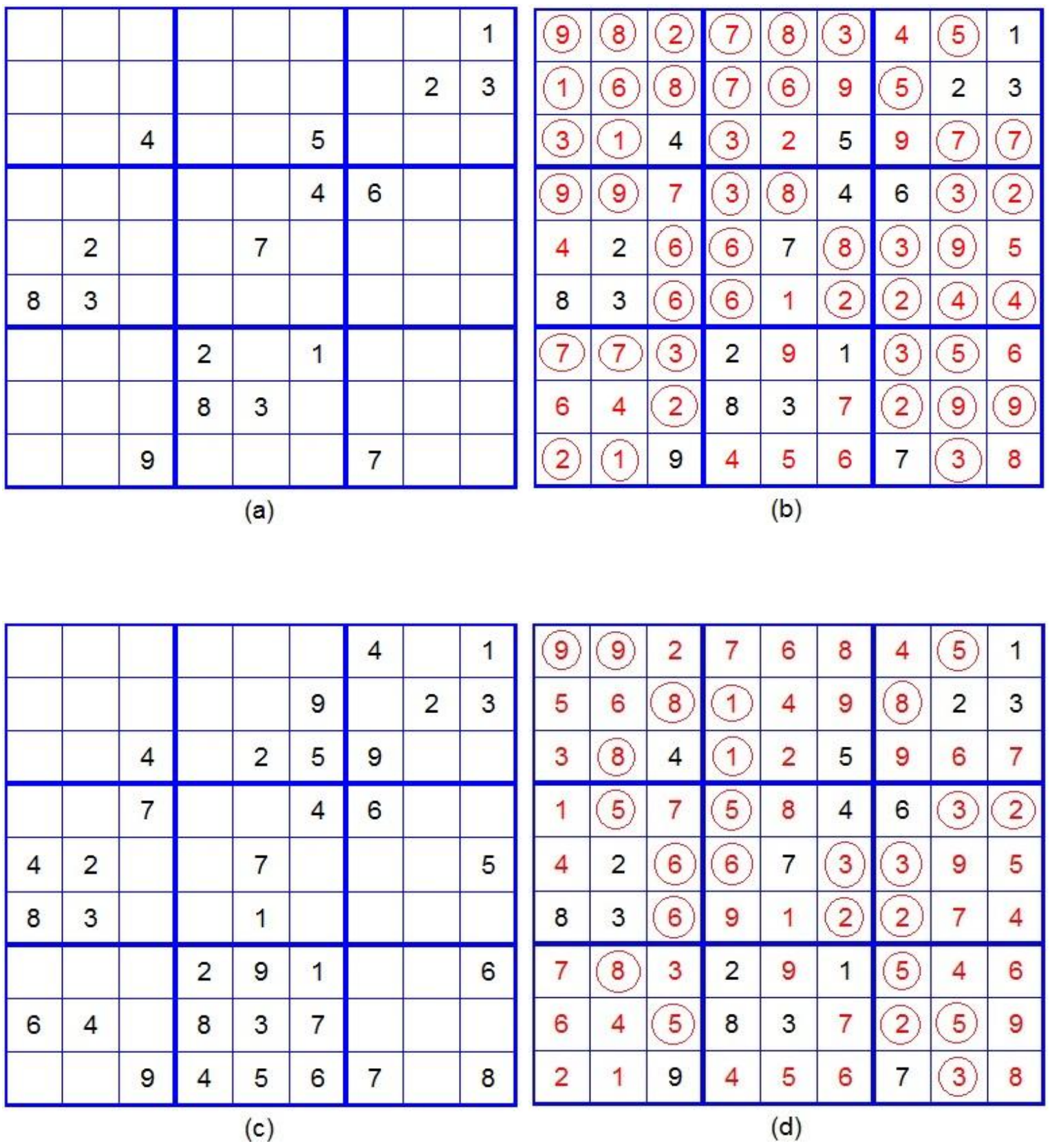

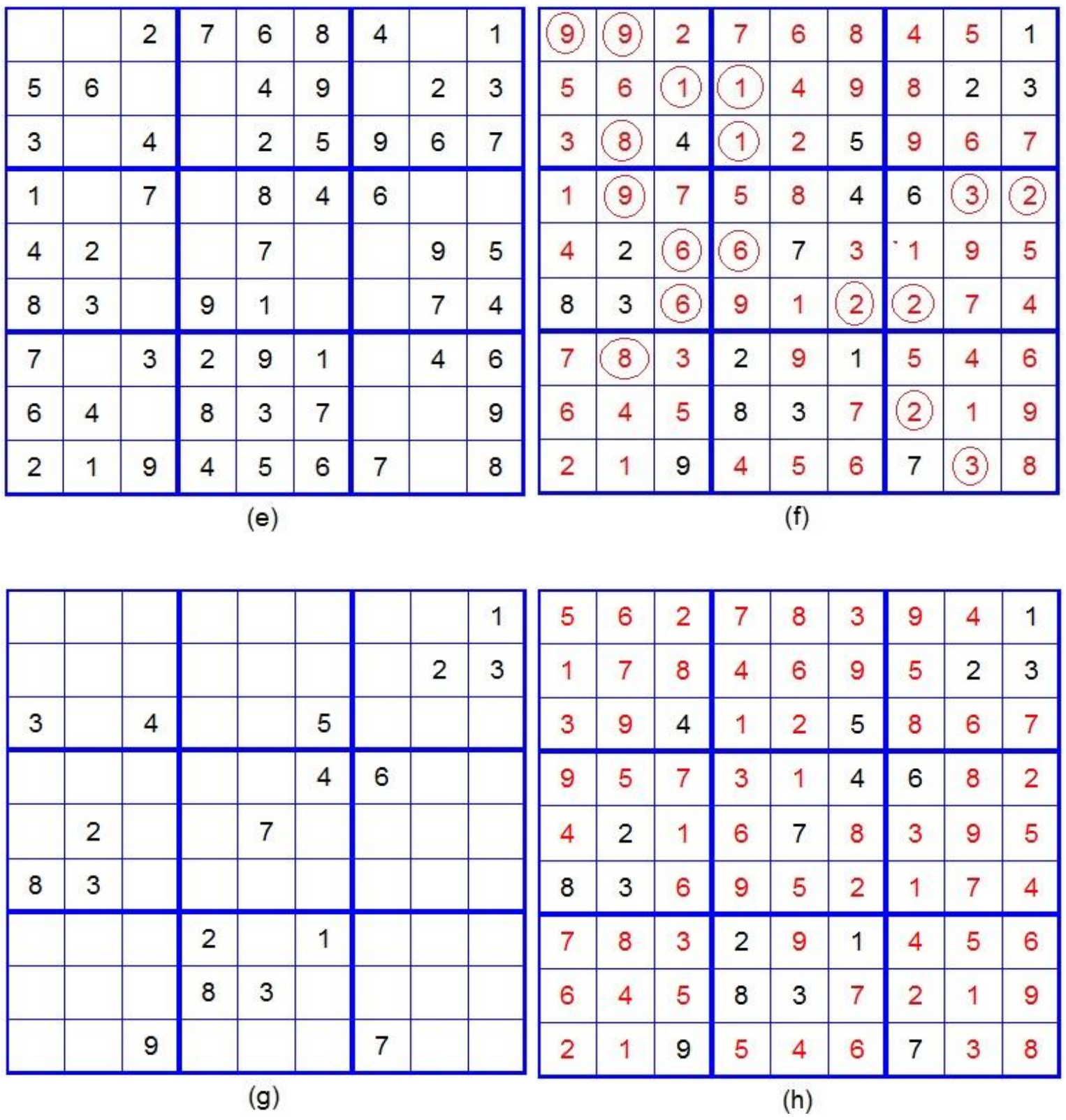

Figure 4. (a) The original Sudoku puzzle; (b) First solved solution, which is wrong; (c) New Sudoku puzzle obtained by deleting the wrong numbers in (b); (d) Second solved solution, which is not correct; (e) New Sudoku puzzle obtained by deleting the wrong numbers in (d); (f) Third solved solution, which is still wrong; (g) New Sudoku puzzle obtained by adding a number to the original Sudoku puzzle; (h) Final correct solution.

In step 3, we pick up a number once, and add it to the original Sudoku puzzle. This complexity is linear, because there are maximum of 81 numbers. If we consider picking up two numbers, then this will reduce to another combination problem, and hence the complexity will increase accordingly. Therefore, we stop our method at step 3.

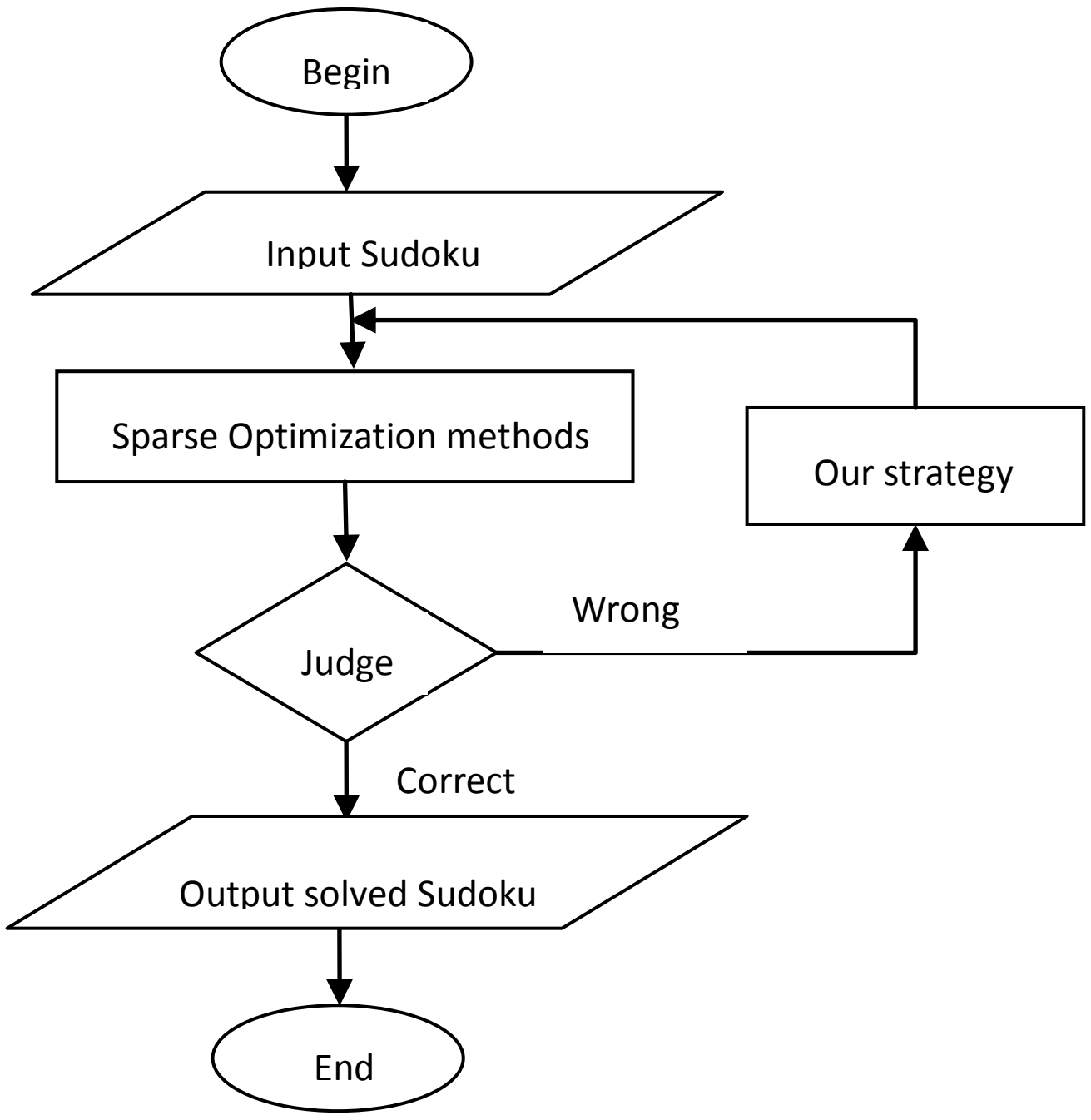

Figure 5. Flowchart of our strategy for solving Sudoku puzzles by sparse optimization methods.

Figure 5 presents a clear overview of our proposed strategy for solving Sudoku based on sparse optimization methods. As a byproduct, we define a new difficulty level of Sudoku as follows:

1) Easy: The given Sudoku is solved directly by the sparse optimization model ($P_1$) or ($WP_1$);

2) Middle: The given Sudoku is solved when Step 1 and Step 2 are used;

3) Hard: The given Sudoku is solved when the Step 3 is used;

4) Devil: The given Sudoku is still not solved using our proposed strategies.

## 4. Numerical experiments

In this section, we test the efficiency of our proposed methods. All the experiments were run on a standard Lenovo laptop with the Intel Core i7-4712MQ CPU 2.3 GHz CPU and 4GB RAM. The software is MATLAB 2013a. We chose a Sudoku dataset a total of 49,151 numbers. This dataset was downloaded from the website of Professor Gordon Royle [10]. All the Sudoku puzzles in the dataset have 17 clues and a unique solution.

Because the solutions of Sudoku puzzles are bounded to be above one, we consider two constraint sets: 1) Nonnegative constraint, $C = \{x|x \geq 0\}$. 2) Bounded constraint, $C = \{x|0 \leq x \leq 1\}$. The results are reported in Table 2 and Table 3. In these tables, the first number in the "First solving" is the total number of successfully solved Sudoku puzzles and the second number is the corresponding percentage. The two numbers in the "Total" column are similarly defined.

Table 2. The successfully solved Sudoku puzzles numbers by sparse optimization model ($P_1$) with and without our proposed strategy

| Methods | Constraint set C | First Solving | Our strategy | | | Total | Time(s) |
|---|---|---|---|---|---|---|---|
| | | | Step 1 | Step 2 | Step 3 | | |
| (LP1) | $\{x \geq 0\}$ | **41722/ 84.8%** | 2246 | 16 | 4713 | **48697/ 99.08%** | 2.2676e+04 |
| | $\{0 \leq x \leq 1\}$ | **41722/ 84.8%** | 2204 | 19 | 4755 | **48700/ 99.08%** | 2.9375e+04 |
| (LP2) | $\{x \geq 0\}$ | **41722/ 84.8%** | 2254 | 22 | 4686 | **48684/ 99.05%** | 5.5224e+03 |
| | $\{0 \leq x \leq 1\}$ | **41722/ 84.8%** | 2180 | 20 | 4747 | **48669/ 99.02%** | 6.7735e+03 |

We can see from Table 2 that the total solved Sudoku puzzles are the same when using the sparse optimization model ($P_1$) with different constraint sets. This only achieves an 84.8% precision, which means that 15.2% of Sudoku puzzles are unsolved. By using our proposed methods, we see that the total number of solved Sudoku puzzle is higher than before, and exceeds 99%. The highest rate is obtained by applying the method (LP1) with the constraint set $0 \leq x \leq 1$.

For the weighted $\ell_1$-norm minimization problem ($WP_1$), we find that the $\varepsilon$ value affects the performances of the corresponding optimization algorithms. For the original ($WP_1$), the $\varepsilon$ value is limited to the [0,1] interval. We then tried large $\varepsilon$ values and the total number of solved Sudoku puzzles in the first step increased accordingly. The maximum number of solved Sudoku puzzles is obtained by applying (Weighted LP2) with the bounded constraint. We can also confirm that the weighted $\ell_1$-norm minimization problem ($WP_1$) outperforms the $\ell_1$-norm minimization problem ($P_1$) for solving Sudoku puzzles.

Table 3. The successfully solved Sudoku puzzles numbers by sparse optimization model ($WP_1$) with and without our proposed strategy

| Methods | Constraint set C | ε | First Solving | Our strategy | | | Total | Time(s) |
|---|---|---|---|---|---|---|---|---|
| | | | | Step 1 | Step 2 | Step 3 | | |
| (Weighted LP1) | $\{x \geq 0\}$ | 0.5 | **45845/ 93.27%** | 80 | 0 | 3006 | **48931/ 99.55%** | 1.3382e+04 |
| | | 1 | **45935/ 93.46%** | 52 | 5 | 2953 | **48945/ 99.58%** | 1.9665e+04 |
| | | 30 | **46028/ 93.65%** | 85 | 0 | 2763 | **48876/ 99.44%** | 2.1967e+04 |
| | $\{0 \leq x \leq 1\}$ | 0.5 | **45683/ 92.94%** | 37 | 1 | 3194 | **48915/ 99.52%** | 1.8959e+04 |
| | | 1 | **45836/ 93.26%** | 105 | 0 | 2975 | **48916/ 99.52%** | 1.9078e+04 |
| | | 30 | **45914/ 93.41%** | 114 | 0 | 2838 | **48866/ 99.42%** | 2.5069e+04 |
| (Weighted LP2) | $\{x \geq 0\}$ | 0.5 | **45858/ 93.30%** | 167 | 0 | 2879 | **48904/ 99.50%** | 2.4295e+04 |
| | | 1 | **45931/ 93.45%** | 110 | 0 | 2837 | **48878/ 99.44%** | 2.7679e+04 |
| | | 30 | **46006/ 93.60%** | 130 | 0 | 2785 | **48921/ 99.53%** | 2.9375e+04 |
| | $\{0 \leq x \leq 1\}$ | 0.5 | **45683/ 92.94%** | 121 | 1 | 3092 | **48897/ 99.48%** | 2.3676e+04 |
| | | 1 | **45839/ 93.26%** | 163 | 0 | 2949 | **48951/ 99.59%** | 2.2516e+04 |
| | | 30 | **45948/ 93.48%** | 114 | 0 | 2893 | **48955/ 99.60%** | 2.2676e+04 |

## 5. Conclusion

The sparse optimization method is a new method for solving Sudoku puzzles. We have proposed an effective strategy for improving this method. This idea is simple and easy to implement without increasing the complexity of the problem. We tested our method on a large dataset of Sudoku puzzles. Numerical results showed that 99%+ of Sudoku puzzles are solved using our method. Although we cannot reach 100% precision, we have enhanced the performance of the original sparse optimization method from 84%+ to 99%+. Furthermore, we presented a new definition of difficulty levels for Sudoku puzzles.

We believe that it may be possible to reach 100% precision by incorporating some logic techniques into sparse optimization methods. However, this exceeds the scope of this paper. We will consider such a method in the future work.

## 6. References

[1] S.S. Skiena, The algorithm design manual, Spinger-Verlag, $2^{nd}$ ed., 2008.

[2] R. Lewis, Metaheuristics can solve Sudoku puzzles, Journal of Heuristics, 13(4), 2007, 387-401.

[3] A. C. Bartlett, T.P. Chartier, A.N. Langville, T.D. Rankin, An integer programming model for the Sudoku problem, The Journal of Online Mathematics and its Applications, 8, 2008, Article ID 1798.
[4] T. K. Moon, J. H. Gunther, J. Kupin, Sinkhorn solves Sudoku, IEEE Transactions on Information Theory, 55( 4) , 2009, 1741–1746.
[5] P. Babu, K. Pelckmans, P. Stoica, J. Li, Linear systems, sparse solutions, and Sudoku, IEEE Signal Processing Letters, 17(1), 2010, 40-42.
[6] Y.D. Zhang, S.H. Wang, Y.K. Huo, L.N. Wu, A novel Sudoku solving methods based on sparse optimization, Journal of Nanjing University of Information Science and Technology: Natural Science Edition, 3(1), 2011, 23-47.
[7] J. Schaad, Modeling the 8-Queens problem and Sudoku using an algorithm based on projections onto nonconvex sets, Master of Science, The University of British Columbia, 2010.
[8] G. McGuire, B. Tugemann, G. Civario, There is no 16-clue Sudoku: solving the Sudoku minimum number of clues problem via hitting set enumeration, Experimental Mathematics, 23(2), 2014, 190-217.
[9] E. Candes, M. Wakin, S. Boyd, Enhancing sparsity by reweighted l1 minimization, Journal of Fourier Analysis Applications, 14(5), 2008, 877-905.
[10] http://staffhome.ecm.uwa.edu.au/~00013890/sudokumin.php.